\documentclass{article}
\usepackage{latexsym}

\usepackage{amsmath,amssymb,amsfonts,array}

\author{Klaus Mohnke\thanks{Universit\"at-Gesamthochschule Siegen,
e-mail:
mohnke@mathematik.uni-siegen.de} }
\title{Legendrian links of topological unknots}
\date{}
\newtheorem{theorem}{Theorem}
\newtheorem{definition}[theorem]{Definition}
\newtheorem{proposition}[theorem]{Proposition}
\newtheorem{lemma}[theorem]{Lemma}
\newtheorem{remark}[theorem]{Remark}

\newtheorem{questions}[theorem]{Questions}
\input{epsf}

\begin{document}
\bibliographystyle{plain}

\maketitle

\begin{abstract}
We use an estimate on the Thurston--Bennequin invariant of a
Legendrian link in terms of its Kauffman--polynomial to show that
links of topological unknots, e.g.~the Borromean rings or the
Whithead link, may not be represented by Legendrian links of
Legendrian unknots.
\end{abstract}

In \cite{Eliashberg:Legendre} Eliashberg classified all Legendre
knots representing the unknot in terms of their Thurston--Bennequin
number, $tb$, and their rotation, $r$. Bennequin's inequality in
these cases reads as
$$
tb + |r|\le -1.
$$
Thus  the Legendre knot given by the wavefront of the 'eye' which
has $tb=-1, r=0$ is referred to as the {\em trivial Legendrian
knot}.

\begin{figure}[h]
\centerline{\epsfbox{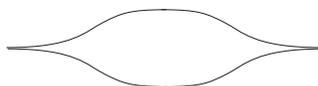}}
\caption{A front projection of the Legendrian unknot}
\label{eye}
\end{figure}

Back then there was no other obvious obstruction  for Legendrian
{\em links} consisting of (topological) unknots and Eliashberg
asked the question: ''Given a link of topological unknots, can it be
realized as a link of [\dots] Legendrian unknots?''

The answer is negative in general and the new obstructions are
given by a sharper inequality on the Thurston--Bennequin number
governed by the Kauffman polynomial $K(x,t)$ which was
found by Lee Rudolph in \cite{Rudolph} (for further reference see
e.g.~\cite{Chmutov/Goryunov} and \cite{Fuks/Tabachnikov}). It simply states
that the Thurston--Bennequin number is not bigger than the
the minimal degree in the variable $x$ of the Kauffman polynomial:
$$
tb\le -\max\mbox{-}\deg_x K.
$$

The contribution of the author is to apply this to links of
topological unknots. We had to be careful because the two groups of
authors \cite{Chmutov/Goryunov,Fuks/Tabachnikov} used different
Kauffman polynomials and thus obtain slightly different
inequalities: here we work with the Dubrovnik--version Chmutov and Goryunov
used.
Let us first recall the definition of the Thurston--Bennequin
number:
\begin{definition}
Let $L$ be an oriented Legendrian link given by a wave front
projection. Then the {\em Thurston--Bennequin number} of $L$,
$tb(L)$, is the number of sideward crossings minus the number of
up-- or downward crossings minus half the number of cusps.

\begin{figure}[h]
\centerline{\epsfbox{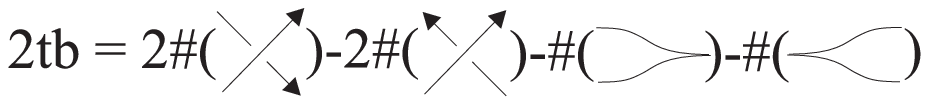}}
\caption{Combinatorial definition of the Thurston--Bennequin number in
terms of the front projection}
\label{tb}
\end{figure}

\end{definition}
From that the following observation is immediate
\begin{lemma}
Let $L=\coprod_iL_i$ be a Legendrian link with pairwise unlinked
components $L_i$. Then the Thurston--Bennequin number of that link
is simply given by the sum of those of the components
$$
tb(L) = \sum_i tb(L_i).
$$
In particular, it does not depend on the orientation of the
components.
\end{lemma}

To investigate Eliashberg's question we took the most simplest
examples we knew: The Borromean rings $B$ and the Whitehead link
$W$. The Kauffman polynomials are given by
\begin{displaymath}\split
K_W(x,y)  &=
yx^5-2x^4-(2y^3+6y)x^3+(-y^4-y^2+6+y^{-2})x^2+(3y^3+9y+2y^{-1})x\\
&\qquad+(y^4+y^2-5-2y^{-2})-(y^3+4y-2y^{-1})x^{-1}+(2+y^{-2})x^{-2}
\endsplit
\end{displaymath}
and
\begin{displaymath}\split
K_B(x,t) & =y^2x^4+(-4y+y^{-3})x^3+(-3y^4-10y^2+3y^{-2})x^2\\
&\qquad+(-2y^5-2y^3+14y+3y^{-1}-3y^{-3})x + (6y^4+18y^2+1-6y^{-2})\\
&\qquad-(-2y^5-2y^3+14y+3y^{-1}-3y^{-3})x^{-1}+(-3y^4-10y^2+3y^{-2})x^{-2}\\
&\qquad-(-4y+y^{-3})x^{-3}+y^2x^4.
\endsplit
\end{displaymath}
From that we  easily deduce the main result of that note
\begin{proposition}
(1) For any Legendrian representaion of the Borromean rings we have
$$
tb\le -4.
$$
(2) For any Legendrian representation of the Whitehaed link we have
$$
tb\le -5.
$$
Thus for both not all components may be Legendrian unknots.
\end{proposition}
\begin{remark}
(1) In the cases of the Borromean rings and the Whitehead link the
inequality is sharp as  the following pictures show. \\
(2) In the case of the Whitehead
link the Kauffman polynomial gives no further obstruction for its
mirror. Indeed it can be represented as a Legendrian link of
Legendrian unknots.

\begin{figure}[h]
\centerline{\epsfbox{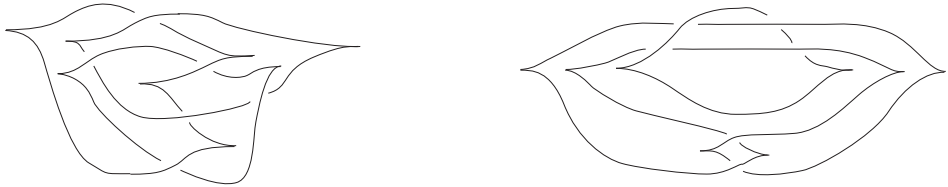}}
\caption{Legendrian Borromean rings with $tb=-4$ and Legendrian Whitehead
mirror consisting of Legendrian unknots}
\label{whiteheadmirror}
\end{figure}

On the other
hand it is possible to give two different Legendrian Whithead links
with $tb=-5$. One consists of components with Thurston--Bennequin
numbers $-4$ and $-1$ the other with $-3$ and $-2$.

\begin{figure}[h]
\centerline{\epsfbox{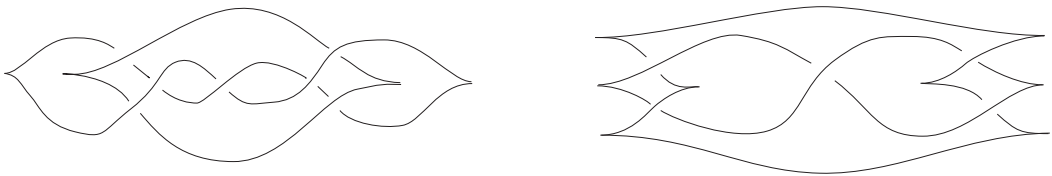}}
\caption{Two Legendrian Whitehead links with $tb=-5$}
\label{whitehead}
\end{figure}

\end{remark}

\begin{questions}
(1) Are there sharper bounds on the Thurston--Bennequin number than
that given in \cite{Chmutov/Goryunov,Fuks/Tabachnikov}?\\ (2) Are
there sharper bounds if one imposes additional conditions?
E.g.~consider Brunnian links, i.e.~links of unknots which fall appart
if one removes one component.\\
(3) Is it at least true that a link of the type described in (2) may be
represented as a Legendrian link consisting of Legendrian unknots
iff the maximal degree in $x$ of its Kauffman polynomial is equal
to the number of components of the link (note that this degree is
never less than the number of components)?\\ (4) Are there further
restrictions for the distribution of the Thurston--Bennequin number
on the components of a Legendrian link?
\end{questions}

{\em Acknowledgements.} I would like to thank Uwe Kaiser and
Alexander Pilz for patiently checking my computations. I apologize
to Lee Rudolph for wrongly stating the reference for the
inequality used in this note in the first version and thank him
for correcting me.

\end{document}